\documentclass[12pt]{article}
\usepackage{mathrsfs}
\usepackage{amssymb}
\usepackage{amsmath,amsthm}

\def\opn#1#2{\def#1{\operatorname{#2}}} 

\newtheorem{Theorem}{Theorem}[section]
\newtheorem{Lemma}[Theorem]{Lemma}
\newtheorem{Corollary}[Theorem]{Corollary}

\newtheorem{Definition}[Theorem]{Definition}

\theoremstyle{definition}
\newtheorem*{pf}{Proof}

\let\epsilon=\varepsilon
\let\kappa=\varkappa
\baselineskip 5pt \lineskip 5pt \numberwithin{equation}{section}
\opn\Aut{Aut} \opn\dist{dist} \opn\mod{mod}
 \opn\ini{in} \opn\inm{inm} \opn\Sym{Sym}
\opn\diag{diag} \opn\Ii{(i)} \opn\Iii{(ii)} \opn\Iif{ if }
\opn\id{id} \opn\Iiii{(iii)} \opn\Iiv{(iv)}
\begin{document}


\title{On the isomorphism classes of Legendre elliptic curves over finite fields
\thanks{Supported by NSF of China (No. 10990011)}}
\author{Rongquan Feng$^{1}$, Hongfeng Wu$^{2}$
\thanks{China Postdoctoral Science Foundation funded project.}
\\ \\
{\small 1~LMAM, School of Mathematical Sciences, Peking
University,}\\{\small Beijing
100871, P.R. China}\\
{\small 2~Academy of Mathematics and Systems Science, Chinese Academy of Sciences,}\\{\small Beijing 100190, P.R. China}\\
{\small\small fengrq@math.pku.edu.cn, whfmath@gmail.com}}

\date{}
\maketitle

\begin{abstract}
In this paper the number of isomorphism classes of Legendre elliptic
curves over finite fields is enumerated.
\end{abstract}


\section{Introduction}
A projective curve is a projective variety of dimension $1$. Let $K$
be an arbitrary field, an elliptic curve $E$ over $K$ is an
absolutely irreducible smooth projective curve of genus $1$ defined
over $K$ with a specified base point $\mathcal{O}$. The elliptic
curve cryptosystem was proposed by Koblitz \cite{Koblitz} and by
Miller \cite{Miller} which relies on the difficulty of discrete
logarithmic problem on the group of rational points on an elliptic
curve.

In order to study the elliptic curve cryptosystem, one need first to
answer how many curves there are up to isomorphism, because two
isomorphic elliptic curves are the same in the point of
cryptographic view. So it is natural to count the isomorphism
classes of some kinds of elliptic curves. Some formulae about
counting the number of the isomorphism classes of general elliptic
curves over a finite field can be found in literatures. For example,
Schoof present the number of isomorphism classes of elliptic curves
over the finite field $\mathbb{F}_q$ in \cite{Schoof}, Menezes
present the number of isomorphism classes of elliptic curves of the
forms $y^2=x^3+ax+b$ over the finite field $\mathbb{F}_q$ in
\cite{Menezes}, and Rezaeian Farashahi and Shparlinski \cite{RFSh}
gave the exact formula for the number of distinct elliptic curves
over a finite field (up to isomorphism over the algebraic closure of
the ground field) in the family of Edwards curves.

In this paper the number of isomorphism classes of Legendre elliptic
curves over the finite field $\mathbb{F}_q$ is enumerated.

\section{Background}
It is well-known that every elliptic curve $E$ over a field $K$ can
be written as a Weierstrass equation
$$E:~Y^2+a_1XY+a_3Y=X^3+a_2X^2+a_4X+a_6$$ with coefficients
$a_1,a_2,a_3,a_4,a_6\in K$. The discriminant $\triangle(E)$ and the
$j$-invariant $j(E)$ of $E$ are defined as
$$\triangle(E)=-b_2^2b_8-8b_4^3-27b_6^2+9b_2b_4b_6$$
and
$$j(E)=(b_2^2-24b_4)^3/\triangle(E),$$ where
\begin{equation*}
\begin{array}{rcl}
     b_2 &=& a_1^2+4a_2,\\
     b_4 &=& 2a_4+a_1a_3, \\
     b_6 &=& a_3^2+4a_6,\\
     b_8 &=& a_1^2a_6-a_1a_3a_4+4a_2a_6+a_2a_3^2-a_4^2.
\end{array}
\end{equation*}

Two projective varieties $V_1$ and $V_2$ are isomorphic if there
exist morphisms $\phi: V_1\rightarrow V_2$ and $\varphi:
V_2\rightarrow V_1$, such that $\varphi\circ \phi$ and
$\phi\circ\varphi$ are the identity maps on $V_1$ and $V_2$
respectively. Two elliptic curves are said to be isomorphic if they
are isomorphic as projective varieties. Let $$E_1:~
Y^2+a_1XY+a_3Y=X^3+a_2X^2+a_4X+a_6$$ and $$E_2~:
Y^2+a_1^{'}XY+a_3^{'}Y=X^3+a_2^{'}X^2+a_4^{'}X+a_6^{'}$$ be two
elliptic curves defined over $K$. It is known \cite{Silveman} that
$E_1$ and $E_2$ are isomorphic over $\overline{K}$, or $E_1$ is
$\overline{K}$-isomorphic to $E_2$, if and only if $j(E_1)=j(E_2)$,
where $\overline{K}$ is the algebraic closure of $K$. However (see
\cite{Silveman}), Let $L$ be an extension of $K$, then $E_1$ and
$E_2$ are isomorphic over $L$ if and only if there exist $u,r,s,t\in
L$ and $u\neq 0$ such that the change of variables
$$(X,Y)\rightarrow (u^2X+r,u^3Y+u^2sX+t)$$ maps the equation of $E_1$ to the equation of $E_2$.
Therefore, $E_1$ and $E_2$ are isomorphic over $L$ if and only if
there exists $u,r,s,t\in L$ and $u\neq 0$ such that
\begin{equation*}\left\{
\begin{array}{rcl}
     ua_1^{'} &=& a_1+2s,\\
     u^2a_2^{'} &=& a_2-sa_1+3r-s^2, \\
     u^3a_3^{'} &=& a_3+ra_1+2t,\\
     u^4a_4^{'} &=& a_4-sa_3+2ra_2-(t+rs)a_1+3r^2-2st,\\
     u^6a_6^{'} &=& a_6+ra_4+r^2a_2+r^3-ta_3-t^2-rta_1.
\end{array}\right.
\end{equation*}
For the simplified Weierstrass equations where
$a_1=a_3=a_1^{'}=a_3^{'}=0$, then $E_1$ is $L$-isomorphic to $E_2$
if and only if there exist $u,r\in L$ and $u\neq 0$ such that
\begin{equation}\label{equiso}
\left\{\begin{array}{rcl}
     u^2a_2^{'} &=& a_2+3r, \\
     u^4a_4^{'} &=& a_4+2ra_2+3r^2,\\
     u^6a_6^{'} &=& a_6+ra_4+r^2a_2+r^3.
\end{array}\right.
\end{equation}
The reader is referred to \cite{Silveman} for more results on the
isomorphism of elliptic curves.

\begin{Definition} The Legendre elliptic curve is one whose Weierstrass equation can be written as
$$E_{\lambda}: y^2=x(x-1)(x-\lambda).$$
\end{Definition}

It is clear that the Legendre elliptic curve $E_{\lambda}$ is
nonsingular for $\lambda\neq 0,1$. The points $\mathcal{O}$,
$(0,0)$, $(1,0)$, and $(\lambda,0)$ are all the $2$-division points,
that is, the points whose double are $\mathcal{O}$'s. The
$j$-invariant of $E_{\lambda}$ is
$j(E_{\lambda})=2^8\frac{(\lambda^2-\lambda+1)^3}{\lambda^2(\lambda-1)^2}$.

\section{Enumeration for Legendre curves}
It is well known \cite{Silveman} that two Legendre curves
$E_{\lambda}: y^2=x(x-1)(x-\lambda)$ and $E_{\mu}:
y^2=x(x-1)(x-\mu)$ are isomorphic over $\overline{\mathbb{F}}_q$ if
and only if they have the same $j$-invariant, or
$$\mu\in\left\{\lambda,\frac{1}{\lambda},1-\lambda,\frac{1}{1-\lambda},\frac{\lambda}{\lambda-1},
\frac{\lambda-1}{\lambda}\right\}.$$ Therefore, the map $\lambda
\mapsto j(L_{\lambda})$ is exactly six-to-one unless when
$\lambda\in\{-1,2,\frac{1}{2}\}$, the map is three-to-one, or when
$\lambda^2-\lambda+1=0$, the map is two-to-one. Note that
$\lambda^2-\lambda+1=0$ has a root in $\mathbb{F}_q$ if and only if
$\mathbb{F}^{*}_q$ has an element of order $3$, which is equivalent
to $q\equiv 1$ or $7~(\text{mod}~12)$. Therefore, we have that the
number of $\overline{\mathbb{F}}_q$-isomorphism classes of Legendre
elliptic curves is $\frac{q-2-3-2}{6}+1+1=\frac{q+5}{6}$ when
$q\equiv
 1,7~(\text{mod}~12)$, and is $\frac{q-2-3}{6}+1=\frac{q+1}{6}$ when $q\equiv 5,11~(\text{mod}~12)$ (see \cite{FengWu}).

In the following, the number of $\mathbb{F}_q$-isomorphism classes
of Legendre elliptic curves will be counted.

\begin{Lemma}\label{isothe0}
Let $L_{a,b}:~y^2=x(x-a)(x-b)$ and $L_{d,e}:~y^2=x(x-d)(x-e)$ be two
elliptic curves defined over $\mathbb{F}_q$, where $ab(a-b)\neq 0$
and $de(d-e)\neq 0$. Then $L_{a,b}$ and $L_{d,e}$ are
$\mathbb{F}_q$-isomorphic if and only if there exists $u\in
\mathbb{F}^{*}_q$ such that the set $\{d,e\}$ is equal to one of
$\{\frac{a}{u^2},\frac{b}{u^2}\}$,
$\{\frac{-a}{u^2},\frac{b-a}{u^2}\}$ and
$\{\frac{a-b}{u^2},\frac{-b}{u^2}\}$.
\end{Lemma}

\begin{pf} From Equation (\ref{equiso}), we know that
$L_{a,b}$ and $L_{d,e}$ are $\mathbb{F}_q$-isomorphic if and only if
there exist $u,r\in \mathbb{F}_q$ and $u\neq 0$ such that
\begin{equation*}\left\{
\begin{array}{rcl}
     -u^2(d+e) &=& -(a+b)+3r, \\
     deu^4 &=& ab-2r(a+b)+3r^2,\\
     0 &=& rab-r^2(a+b)+r^3.
\end{array}\right.
\end{equation*}
Thus $r=0$ or $r=a$ or $r=b$. Therefore, $L_{a,b}$ and $L_{d,e}$ are
$\mathbb{F}_q$-isomorphic if and only if there exists $u\in
\mathbb{F}^{*}_q$ such that
\begin{equation*}
\left\{\begin{array}{rcl}
  (d+e)u^2&=&a+b,\\
  deu^4&=&ab,
\end{array}
\right. \text{~or~} \left\{\begin{array}{rcl}
  (d+e)u^2&=&b-2a,\\
  deu^4&=&a^2-ab,
\end{array}
\right. \text{~or~} \left\{\begin{array}{rcl}
  (d+e)u^2&=&a-2b,\\
  deu^4&=&b^2-ab.
\end{array}
\right.
\end{equation*}
Solving these 3 equations, we get that the set $\{d,e\}$ is
$\{\frac{a}{u^2}, \frac{b}{u^2}\}$, or $\{\frac{-a}{u^2},
\frac{b-a}{u^2}\}$ or $\{\frac{a-b}{u^2}, \frac{-b}{u^2}\}$. \qed
\end{pf}

\begin{Corollary}\label{isothe1}
Let $E_{\lambda}:~y^2=x(x-1)(x-\lambda)$ and
$E_{\mu}:~y^2=x(x-1)(x-\mu)$ be two Legendre elliptic curves defined
over $\mathbb{F}_q$. Then $E_{\lambda}$ and $E_{\mu}$ are
$\mathbb{F}_q$-isomorphic if and only if there exists $u\in
\mathbb{F}^{*}_q$ such that the set $\{1,\mu\}$ is equal to one of
$\{\frac{1}{u^2},\frac{\lambda}{u^2}\}$,
$\{\frac{-1}{u^2},\frac{\lambda-1}{u^2}\}$ and
$\{\frac{1-\lambda}{u^2},\frac{-\lambda}{u^2}\}$.
\end{Corollary}

Now let $E_{\lambda}:~y^2=x(x-1)(x-\lambda)$ be a Legendre elliptic
curve with $\lambda(\lambda-1)\neq 0$. We want to determine how many
Legendre elliptic curves which are $\mathbb{F}_q$-isomorphic to
$E_{\lambda}$. Let $E_{\mu}$ be one of them. Then
$j(E_{\mu})=j(E_{\lambda})$. Thus $\mu\in\{\lambda, 1-\lambda,
\frac{1}{\lambda}, \frac{\lambda-1}{\lambda}, \frac{1}{1-\lambda},
\frac{\lambda}{\lambda-1}\}$. For the different choices of $\mu$,
Table 1 gives what the $u^2$ and what the set $\{1,\mu\}$ can be
compared with the result in Corollary \ref{isothe1}.

\begin{table}
\caption{What the $u^2$ and what the set $\{1,\mu\}$ can be}
\begin{center}
\begin{tabular}{|c|c|c|}
  \hline
  $\mu=$ & $u^2=$ & $\{1,\mu\}=$ \\ \hline
  $\lambda$ & $1$ & $\{\frac{1}{u^2},\frac{\lambda}{u^2}\}$\\ \hline
  $1-\lambda$ & $-1$ & $\{\frac{-1}{u^2},\frac{\lambda-1}{u^2}\}$\\ \hline
  $\frac{1}{\lambda}$ & $\lambda$ & $\{\frac{1}{u^2},\frac{\lambda}{u^2}\}$\\ \hline
  $\frac{\lambda-1}{\lambda}$ & $-\lambda$ & $\{\frac{1-\lambda}{u^2},\frac{-\lambda}{u^2}\}$\\ \hline
  $\frac{1}{1-\lambda}$ & $\lambda-1$ & $\{\frac{-1}{u^2},\frac{\lambda-1}{u^2}\}$\\ \hline
  $\frac{\lambda}{\lambda-1}\}$ & $1-\lambda$ & $\{\frac{1-\lambda}{u^2},\frac{-\lambda}{u^2}\}$\\ \hline
\end{tabular}
\end{center}
\end{table}

Throughout the paper, the characteristic of the finite field
$\mathbb{F}_q$ is assumed to be greater than 3. We distinguish the
enumeration
into the following two cases.\\

\noindent {\bf Case 1: $q\equiv 3~(\text{mod}~4)$}\\

In this case, $-1$ is not a square in $\mathbb{F}_q$. When
$j(E_{\lambda})=0$, then $\lambda^2-\lambda+1=0$. The equation
$\lambda^2-\lambda+1=0$ has roots in $\mathbb{F}_q$ if and only if
$q\equiv 7~(\text{mod}~12)$. The roots are
$\varepsilon_1=\frac{1+\sqrt{-3}}{2}$ and
$\varepsilon_2=\frac{1-\sqrt{-3}}{2}$. By Corollary \ref{isothe1},
it is easy to check that $E_{\varepsilon_1}$ is not
$\mathbb{F}_q$-isomorphic to $E_{\varepsilon_2}$. When
$j(E_{\lambda})=1728$, we have $\lambda\in\{-1,2,\frac{1}{2}\}$.
Since $\{1,2\}=\{\frac{1-(-1)}{1},\frac{-(-1)}{1}\}$, $E_{-1}$ is
$\mathbb{F}_q$-isomorphic to $E_{2}$. Furthermore, we know that
exactly one of $-2$ and $2$ is a square. When $2$ is a square in
$\mathbb{F}_q$, then
$\{1,\frac{1}{2}\}=\{\frac{1}{(\sqrt{2})^2},\frac{2}{(\sqrt{2})^2}\}$.
When $-2$ is a square in $\mathbb{F}_q$, then
$\{1,\frac{1}{2}\}=\{\frac{1-2}{(\sqrt{-2})^2},\frac{-2}{(\sqrt{-2})^2}\}$.
Hence $E_{\frac{1}{2}}$ is $\mathbb{F}_q$-isomorphic to $E_{2}$.
That means $E_{-1}$, $E_{2}$ and $E_{\frac{1}{2}}$ are
$\mathbb{F}_q$-isomorphic to each other. Now assume that
$j(E_{\lambda})\neq 0,1728$, then the $6$ possible values of $\mu$
are different to each other. But now exactly one of $1$ and $-1$,
one of $\lambda$ and $-\lambda$, one of $\lambda-1$ and $1-\lambda$
is a square in $\mathbb{F}_q$. So there are exactly $3$ values of
$\mu$ such that $E_{\mu}$ is $\mathbb{F}_q$-isomorphic to
$E_{\lambda}$. Therefore the number of $\mathbb{F}_q$-isomorphism
classes of Legendre curves is
\begin{equation*}
\left\{
\begin{array}{ll}
 \dfrac{q-2-3-2}{3}+1+2=\dfrac{q+2}{3},~~&\text{if~~} q\equiv 7~(\text{mod}~12),\\
 [2ex]
 \dfrac{q-2-3}{3}+1=\dfrac{q-2}{3},~~&\text{if~~} q\equiv 11~(\text{mod}~12).
\end{array}
\right.
\end{equation*}
This proves the following theorem.
\begin{Theorem}\label{lejco1}
Suppose $\mathbb{F}_q$ is a finite field with $char(\mathbb{F}_q)>3$
and $q\equiv 3~(\text{mod}~4)$. Let $N_q$ be the number of
$\mathbb{F}_q$-isomorphism classes of Legendre curves $E_{\lambda}:
y^2=x(x-1)(x-\lambda)$ defined over $\mathbb{F}_q$ with
$\lambda(\lambda-1)\neq 0$. Then
\begin{equation*}
N_q=\left\{
\begin{array}{ll}
 \dfrac{q+2}{3},~~&\text{if~~} q\equiv 7~(\text{mod}~12),\\[2ex]
 \dfrac{q-2}{3},~~&\text{if~~} q\equiv 11~(\text{mod}~12).
\end{array}
\right.
\end{equation*}
\end{Theorem}

\noindent {\bf Case 2: $q\equiv 1~(\text{mod}~4)$}\\

For the finite field $\mathbb{F}_q$, the Jacobi symbol
$\left(\dfrac{a}{q}\right)$ is defined as follows for all
$a\in\mathbb{F}_q$:
\begin{equation*}
\left(\dfrac{a}{q}\right)=\left\{
\begin{array}{rcl}
 1,~~~&&\text{if~~} a\neq 0 \text{~and~} a \text{~is a square in~} \mathbb{F}_q,\\
 -1,~~~&&\text{if~~} a \text{~is not a square in~} \mathbb{F}_q,\\
 0,~~~&&\text{if~~} a=0.
\end{array}
\right.
\end{equation*}
The following lemma can be got easily.

\begin{Lemma}\label{residu1}
Suppose that $\mathbb{F}_q$ is a finite field with
$char(\mathbb{F}_q)>3$ and $q\equiv 1~(\text{mod}~4)$. Let $N(s,t)$
be the number of $a\in\mathbb{F}_q$ with
$\left(\frac{a}{q}\right)=s$ and $\left(\frac{1-a}{q}\right)=t$.
Then $N(1,1)=\frac{q-5}{4}$,
$N(1,-1)=N(-1,1)=N(-1,-1)=\frac{q-1}{4}$.
\end{Lemma}

Now we divide the Legendre elliptic curves
$E_{\lambda}:y^2=x(x-1)(x-\lambda)$ with $\lambda\neq 0,1$, into the
following 4 disjoint sets $H_1$, $H_2$, $H_3$ and $H_4$, where
\begin{equation*}
\begin{array}{rcl}
 H_1&=&\left\{y^2=x(x-1)(x-b)|\left(\frac{b}{q}\right)=\left(\frac{1-b}{q}\right)=1\right\},\\
 [1ex]
 H_2&=&\left\{y^2=x(x-1)(x-b)|\left(\frac{b}{q}\right)=1, \left(\frac{1-b}{q}\right)=-1\right\},\\
 [1ex]
 H_3&=&\left\{y^2=x(x-1)(x-b)|\left(\frac{b}{q}\right)=-1, \left(\frac{1-b}{q}\right)=1\right\},\\
 [1ex]
H_4&=&\left\{y^2=x(x-1)(x-b)|\left(\frac{b}{q}\right)=-1,
\left(\frac{1-b}{q}\right)=-1\right\}.
\end{array}
\end{equation*}
From Lemma \ref{residu1},  we get that $|H_1|=\frac{q-5}{4}$ and
$|H_2|=|H_3|=|H_4|=\frac{q-1}{4}$.

Note that $-1$ is a square in this case. As in Case 1,
$j(E_{\lambda})=0$ if and only if $\lambda^2-\lambda+1=0$. The
equation $\lambda^2-\lambda+1=0$ has roots in $\mathbb{F}_q$ if and
only if $q\equiv 1,13~(\text{mod}~24)$. The roots are
$\varepsilon_1=\frac{1+\sqrt{-3}}{2}=\left(\frac{\sqrt{3}+\sqrt{-1}}{2}\right)^2$
and
$\varepsilon_2=\frac{1-\sqrt{-3}}{2}=\left(\frac{\sqrt{3}-\sqrt{-1}}{2}\right)^2=1-\varepsilon_1$.
Thus $E_{\varepsilon_1}, E_{\varepsilon_2}\in H_1$ and by Corollary
\ref{isothe1}, $E_{\varepsilon_1}$ is $\mathbb{F}_q$-isomorphic to
$E_{\varepsilon_2}$ since
$\left\{1,\varepsilon_2\right\}=\left\{\frac{-1}{(\sqrt{-1})^2},\frac{\varepsilon_1-1}{(\sqrt{-1})^2}\right\}$.
When $j(E_{\lambda})=1728$, we have
$\lambda\in\{-1,2,\frac{1}{2}\}$. Now $2$ is a square if and only if
$q\equiv 1,17~(\text{mod}~24)$. So when $q\equiv
1,17~(\text{mod}~24)$, $E_{-1}, E_{2}, E_{\frac{1}{2}}\in H_1$ and
they are $\mathbb{F}_q$-isomorphic to each other as in Case 1. When
$q\equiv 5,13~(\text{mod}~24)$, $E_{-1}\in H_2$, $E_{2}\in H_3$,
while $E_{\frac{1}{2}}\in H_4$. It can be checked easily that
$E_{-1}$ and $E_{2}$ are $\mathbb{F}_q$-isomorphic to each other but
$E_{\frac{1}{2}}$ is not $\mathbb{F}_q$-isomorphic to them.

Now let $j(E_{\lambda})\neq 0,1728$. Then $\lambda, 1-\lambda,
\frac{1}{\lambda}, \frac{\lambda-1}{\lambda}, \frac{1}{1-\lambda},
\frac{\lambda}{\lambda-1}$, the 6 possible values for $\mu$, are
distinct with each other. If $E_{\lambda}\in H_1$, then all the 6
values $1, -1, \lambda,-\lambda, 1-\lambda$ and $\lambda-1$ are
squares in $\mathbb{F}_q$. Thus for any $\mu\in\{\lambda, 1-\lambda,
\frac{1}{\lambda}, \frac{\lambda-1}{\lambda}, \frac{1}{1-\lambda},
\frac{\lambda}{\lambda-1}\}$, $E_{\mu}$ is $\mathbb{F}_q$-isomorphic
to $E_{\lambda}$ and $E_\mu\in H_1$. If $E_{\lambda}\in H_4$, then
$E_{\mu}$ is $\mathbb{F}_q$-isomorphism to $E_{\lambda}$ only for
$\mu\in\{\lambda,1-\lambda\}$ since only $1$ and $-1$ are squares
among the six numbers in the second column of Table 1. Now
$E_{1-\lambda}\in H_4$ too.

If $E_{\lambda}\in H_2$, then $1,-1, \lambda$ and $-\lambda$ are
squares while $1-\lambda$ and $\lambda-1$ are non-squares. Thus
$E_{\mu}$ is $\mathbb{F}_q$-isomorphic to $E_{\lambda}$ for
$\mu\in\{\lambda, 1-\lambda, \frac{1}{\lambda},
\frac{\lambda-1}{\lambda}\}$. Note that now
$E_{\frac{1}{\lambda}}\in H_2$ but $E_{1-\lambda},
E_{\frac{\lambda-1}{\lambda}}\in H_3$. Similarly, if $E_{\lambda}\in
H_3$, then $E_{\mu}$ is $\mathbb{F}_q$-isomorphic to $E_{\lambda}$
for $\mu\in\{\lambda, 1-\lambda, \frac{1}{1-\lambda},
\frac{\lambda}{\lambda-1}\}$ and $E_{\frac{\lambda-1}{\lambda}}\in
H_3$ but $E_{1-\lambda}, E_{\frac{1}{1-\lambda}}\in H_2$.

Therefore, let $N_{q,H_1}$, $N_{q,H_2\cup H_3}$, and $N_{q,H_4}$ be
the number of $\mathbb{F}_q$-isomorphism classes of Legendre
elliptic curves $E_{\lambda}\in H_1$, $H_2\cup H_3$ and $H_4$
respectively. Then we have
\begin{equation*}
N_{q,H_1}=\left\{
\begin{array}{ll}
 \left(\dfrac{q-5}{4}-3-2\right)/6+1+1=\dfrac{q+23}{24},~~&\text{if~~} q\equiv 1~(\text{mod}~24),\\
 [2ex]
 \left(\dfrac{q-5}{4}\right)/6=\dfrac{q-5}{24},~~&\text{if~~} q\equiv 5~(\text{mod}~24),\\
 [2ex]
  \left(\dfrac{q-5}{4}-2\right)/6+1=\dfrac{q+11}{24},~~&\text{if~~} q\equiv
 13~(\text{mod}~24),\\
 [2ex]
 \left(\dfrac{q-5}{4}-3\right)/6+1=\dfrac{q+7}{24},~~&\text{if~~} q\equiv 17~(\text{mod}~24),
\end{array}
\right.
\end{equation*}
\begin{equation*}
N_{q,H_2\cup H_3}=\left\{
\begin{array}{ll}
 \left(\dfrac{q-1}{2}\right)/4=\dfrac{q-1}{8},~~&\text{if~~} q\equiv 1,17~(\text{mod}~24),\\
 [2ex]
 \left(\dfrac{q-1}{2}-2\right)/4+1=\dfrac{q+3}{8},~~&\text{if~~} q\equiv
 5,13~(\text{mod}~24),
\end{array}
\right.
\end{equation*}
and
\begin{equation*}
N_{q,H_4}=\left\{
\begin{array}{ll}
 \left(\dfrac{q-1}{4}\right)/2=\dfrac{q-1}{8},~~&\text{if~~} q\equiv
 1,17~(\text{mod}~24),\\[2ex]
 \left(\dfrac{q-1}{4}-1\right)/2+1=\dfrac{q+3}{8},~~&\text{if~~} q\equiv 5,13~(\text{mod}~24).
\end{array}
\right.
\end{equation*}


We know from above analysis that the Legendre curves from the 3
distinct sets $H_1$, $H_2\cup H_3$ and $H_4$ can not be
$\mathbb{F}_q$-isomorphic to each other. Summing up the above
numbers, we have the following theorem:
\begin{Theorem}\label{lejco2}
Suppose $\mathbb{F}_q$ is a finite field with $char(\mathbb{F}_q)>3$
and $q\equiv 1~(\text{mod}~4)$. Let $N_q$ be the number of
$\mathbb{F}_q$-isomorphism classes of Legendre curves $E_{\lambda}:
y^2=x(x-1)(x-\lambda)$ defined over $\mathbb{F}_q$ with
$\lambda(\lambda-1)\neq 0$. Then
\begin{equation*}
N_q=\left\{
\begin{array}{ll}
 \dfrac{7q+17}{24},~~&\text{if~~} q\equiv 1~(\text{mod}~24),\\[2ex]
 \dfrac{7q+13}{24},~~&\text{if~~} q\equiv 5~(\text{mod}~24),\\[2ex]
 \dfrac{7q+29}{24},~~&\text{if~~} q\equiv 13~(\text{mod}~24),\\[2ex]
 \dfrac{7q+1}{24},~~&\text{if~~} q\equiv 17~(\text{mod}~24).
\end{array}
\right.
\end{equation*}
\end{Theorem}

Combining Theorems \ref{lejco1} and \ref{lejco2} together, we have
the following enumeration result.
\begin{Theorem}
Suppose $\mathbb{F}_q$ is the finite field with $q$ elements and
$char(\mathbb{F}_q)>3$. Let $N_q$ be the number of
$\mathbb{F}_q$-isomorphism classes of Legendre curves $E_{\lambda}:
y^2=x(x-1)(x-\lambda)$ defined over $\mathbb{F}_q$ with
$\lambda(\lambda-1)\neq 0$. Then
\begin{equation*}
N_q=\left\{
\begin{array}{ll}
\dfrac{7q+17}{24},~~&\text{if~~} q\equiv 1~(\text{mod}~24),\\[2ex]
 \dfrac{7q+13}{24},~~&\text{if~~} q\equiv 5~(\text{mod}~24),\\[2ex]
 \dfrac{q+2}{3},~~&\text{if~~} q\equiv 7,19~(\text{mod}~24),\\[2ex]
 \dfrac{q-2}{3},~~&\text{if~~} q\equiv 11,23~(\text{mod}~24),\\[2ex]
  \dfrac{7q+29}{24},~~&\text{if~~} q\equiv 13~(\text{mod}~24),\\[2ex]
 \dfrac{7q+1}{24},~~&\text{if~~} q\equiv 17~(\text{mod}~24).
\end{array}
\right.
\end{equation*}
\end{Theorem}


\begin{thebibliography}{90}
\bibitem{FengWu} R. Feng and H. Wu, Number of general Jacobi quartic curves over finite
fields, preprint.

\bibitem{Dale} D. Husem\"{o}ller, Elliptic Curves, Graduate Texts in Mathematics 111, second edition,
Springer-Verlag, New York, 2004.

\bibitem{Koblitz} N. Koblitz, ¡®Elliptic curve cryptosystems¡¯, Math. Comp., 48(177),
(1987), 203¨C209.

\bibitem{Menezes} A.J. Menezes, Elliptic Curve Public Key Cryptosystems,
Kluwer Academic Publishers, 1993.

\bibitem{Miller} V.S. Miller, ¡®Use of elliptic curves in cryptography¡¯, Advances
in Cryptology ¨C Crypto 1985, Lecture Notes in Comp. Sci., vol. 218,
Springer-Verlag, 1986, 417¨C426.

\bibitem{RFSh} R. Rezaeian Farashahi and I. E. Shparlinski. ¡®On the number of
distinct elliptic curves in some families¡¯, Designs, Codes and
Cryptography, 83-99, Vol.54, No.1, 2010.

\bibitem{Schoof} R. Schoof, Nonsigular plane cubic curves over
finite field, J.\ Combine, Theory Ser.\ A 46(1987), 183-211.

\bibitem{Silveman} J.H. Silverman, The Arithmetic of Elliptic Curves, Graduate Texts in Mathematics
106, Springer-Verlag, New York, 1986.

\end{thebibliography}
\end{document}